\documentclass[11pt]{article}

\usepackage{amssymb, amsmath}
\usepackage{float,epsfig}
\usepackage{epsfig,amsfonts,amsbsy,amssymb,amsthm,here}
\usepackage[usenames,dvipsnames]{color}
\usepackage{graphics}
\usepackage{graphicx}
\usepackage{epsfig}
\usepackage{url}
\usepackage{arydshln,leftidx,mathtools}
\usepackage{rotating}
\usepackage{lscape}
\usepackage{setspace}
\usepackage{algpseudocode}
\usepackage{algorithm}
\usepackage[all]{xy}

\newtheorem{theorem}{\bf Theorem}[section]
\newtheorem{proposition}[theorem]{\bf Proposition}
\newtheorem{definition}[theorem]{\bf Definition}
\newtheorem{corollary}[theorem]{\bf Corollary}

\newtheorem{example}[theorem]{\bf Example}

\newtheorem{remark}[theorem]{\bf Remark}
\newtheorem{lemma}[theorem]{\bf Lemma}

\newcommand{\ba}{\begin{array}}
\newcommand{\ea}{\end{array}}

\newcommand{\beano}{\begin{eqnarray*}}
\newcommand{\eeano}{\end{eqnarray*}}

\def\bmatrix#1{\left[ \begin{matrix} #1 \end{matrix} \right]}

\def \R{{\mathbb R}}
\def \C{{\mathbb C}}

\def \bmatrix#1{\left[ \begin{matrix} #1 \end{matrix} \right]}
\def \H{{\mathbb H}}

\def \Delta{\triangle}

\def \rank{\mathrm{rank}}
\def \pf{{\bf Proof: }}

\def \rank{\mathrm{rank}}

\def \pf{{\bf Proof: }}

\def \x{{\mathbf{x}}}

\def \0{{\mathbf{0}}}

\def \rq{{\mathbf{r}}}
\def \s{{\mathbf{s}}}
\def \t{{\mathbf{t}}}
\def \u{{\mathbf{u}}}

\title{\textbf{A Complete Characterization of Determinantal Quadratic Polynomials}} 
\begin{document}
\author{Papri Dey\thanks{Email:papridey@ee.iitb.ac.in, yedirpap@gmail.com} \,\, and \,\,  Harish K. Pillai\thanks{Email: hp@ee.iitb.ac.in} \\ Department of Electrical Engineering\\ IIT Bombay, India}
\textheight 8.9in
 \textwidth 6.1in
\topmargin -0.5in \evensidemargin .2in \oddsidemargin .3in

\date{}
\maketitle
\begin{abstract}
The problem of expressing a multivariate polynomial as the determinant of a monic (definite) symmetric or Hermitian linear matrix polynomial (LMP) has drawn a huge amount of attention due to its connection with optimization problems. In this paper we provide a necessary and sufficient condition for the existence of \textit{monic Hermitian determinantal representation} as well as \textit{monic symmetric determinantal representation} of size $2$ for a given quadratic polynomial. Further we propose a method to construct such a monic determinantal representtaion (MDR) of size $2$ if it exists. It is known that a quadratic polynomial $f(\x)=\x^{T}A\x+b^{T}\x+1$ has a symmetric MDR of size $n+1$ if $A$ is \textit{negative semidefinite}. We prove that if a quadratic polynomial $f(\x)$ with $A$ which is not negative semidefinite has an MDR of size greater than $2$, then it has an MDR of size $2$ too. We also characterize quadratic polynomials which exhibit diagonal MDRs. 
\end{abstract}
\section{Introduction}
This paper deals with characterization of quadratic real multivariate polynomials which admit monic Hermitian (symmetric) determinantal representations, that is polynomials which can be written as the determinant of a monic linear matrix polynomial (LMP) whose coefficient matrices are Hermitian (symmetric) and the constant coefficient matrix is the identity matrix. Note that the coefficient matrices of the LMP could be of any order greater than two. In particular, in this paper, we focus on the existence and determination of a monic LMP whose coefficient matrices are Hermitian (symmetric) and of order $2$ for a given quadratic polynomial. Besides, we identify the class of quadratic polynomials for which an MDR of size greater than $2$ ensures the existence of an MDR of size $2$ respectively.

Determinantal representations of polynomials have generated a lot of interest due to its connection with the problem of determining (definite) LMI representable sets, also known as spectrahedra \cite{pablo2} which play a crucial role in optimization problems. Indeed, if the feasible set of an optimization problem is a definite LMI representable set, the problem can be transformed into a semidefinite programming (SDP) problem which in turn can be efficiently solved by SDP solvers. It is important to recall that any polynomial can be expressed as the determinant of a symmetric LMP \cite{Quarez}.
It is also known that the algebraic interior defined by a real zero (RZ) quadratic polynomial is always a spectrahedron, since a Hermitian determinantal representation can be obtained for higher powers of RZ quadratic polynomials using Clifford algebra \cite{Tim} and sum of squares (SOS) decomposition of a parametrized Hermite matrix \cite{Timhermite}. To the best of authors' knowledge, characterization of quadratic polynomials that have an MDR of size $2$ has not been done before.

In this paper, we provide a necessary and sufficient condition for the existence of MDRs of size $2$ for any quadratic polynomial. We also propose a constructive method to determine these MDRs. We show that for a certain sub-class of quadratic polynomials that have MDRs of size $2$, there are precisely two non-equivalent classes of MDRs, whereas for all other quadratic polynomials that have MDRs of size $2$, all MDRs are unitarily equivalent. Recall
that a quadratic polynomial $f(\x) = \x^{T}A\x + b^{T} \x + 1$ admits a symmetric MDR of size $n+1$ if the corresponding matrix $A$ is negative semidefinite \cite{Ramana1}, \cite{Nemirovski1}. This need not imply the existence of MDR of size $2$ for the same polynomial. Therefore, it is natural to
ask whether quadratic polynomials which have an MDR of size $2$ can be characterized. We further characterize all quadratic polynomials that have an MDR of any size [Sec \ref{seccompmsdr}]. We show that if a quadratic polynomial $f(\x) = \x^{T}A\x + b^{T} \x + 1$ has an MDR, then either $A$ is negative semidefinite or $f(\x)$ admits an MDR of size $2$. In other words, if $f(\x)$ has an MDR, but $A$ is not negative semidefinite, then $f(\x)$ has an MDR of size $2$.
\section{Preliminaries}
We begin with the concept of definite LMI representable sets and its relation to monic determinantal representations. A set $S \subseteq \R^{n}$ is said to be \textit{LMI representable} if
\begin{equation} \label{lmiset}
S= \{ \x \in \R^{n} : A_{0} +x_{1}A_{1} +x_{2}A_{2} + \dots + x_{n}A_{n} \succeq 0 \}
\end{equation}
for some real symmetric matrices $A_{i}, i=0,\dots,n$ and $\x=(x_{1}, \dots, x_{n})^{T}$. If $A_{0}\succ 0$, the set $S$ is called a \textit{definite LMI representable} set whereas if $A_{0}=I$, $S$ is known to be a \textit{monic LMI representable} set. By $A \succ 0 (\succeq 0)$ we mean that the matrix $A$ is positive (semi)-definite. A \textit{spectrahedron} is another name used for an LMI representable set. It is evident that a spectrahedron is both convex and \textit{basic closed semialgebraic} set. Moreover, if a spectrahedron has a nonempty interior, it is a definite LMI representable set [\cite{Ramana1},section $1.4$], \cite{Tim} -- without loss of generality, the origin may be considered as an interior point of the set.  It is also known that a definite LMI representable set is always monic LMI representable \cite{Helton}.

A polynomial $f(\x) \in \R[\x]$ is said to have a determinantal representation if $f(\x)$ is the determinant of a linear matrix polynomial (LMP) i.e.,
\begin{equation} \label{detreppoly}
f(\x) = \det(A_{0}+x_{1}A_{1} +x_{2}A_{2} + \dots + x_{n}A_{n}), \mbox{ where } A_{i} \in \H^{k \times k}(\C) \ \mbox{for some} \ k.
\end{equation}
If the matrices $A_{i} \in S\R^{k \times k}$ (symmetric matrices of order $k$), then the polynomial has a symmetric determinantal representation. Note that as $f(\x) \in \R[\x]$, the matrices $A_{i}$ could have been Hermitian matrices too. Therefore, if the matrices $A_{i} \in \H^{k \times k}(\C)$ (Hermitian matrices of order $k$), then the polynomial is said to have a Hermitian determinantal representation. When the matrices $A_{i}$s are of size $k$, we call $k$ the size of the determinantal representation.

The determinantal representation is definite if $A_{0} \succ 0$. Further if $A_{0} = I_{k}$, the identity
matrix of order $k$, then we have a monic determinantal representation (MDR). Throughout the paper, we are interested in monic determinantal representations of polynomials using either symmetric or Hermitian matrices. Therefore, we use the acronyms MSDR and MHDR for monic symmetric determinantal representation and monic Hermitian determinantal representation, respectively. If all $A_{i} \in S\R^{k \times k}$ are diagonal, then the polynomial is said to have
a diagonal determinantal representation. It is obvious that a polynomial which admits a definite determinantal representation can be scaled in order to admit a monic determinantal representation. Without loss of generality, throughout this paper we consider only problems dealing with monic determinantal representations and hence we consider only those polynomials whose constant coefficient is $1$, unless stated otherwise.

If a polynomial $f(\x)$ also admits an MDR, i.e., $f(\x)=\det(I+x_{1}A_{1} +x_{2}A_{2} + \dots + x_{n}A_{n})$, then $f(\x)>0$ when $\x\in\mbox{Int}(S)$ where the spectrahedron $S=\{\x : I+x_{1}A_{1} +x_{2}A_{2} + \dots + x_{n}A_{n}\succeq 0\}$ and $f(\x)=0$ when $\x\in\partial S$. On the other hand, given $f(\x) \in \R[\x]$, a closed subset $C_{f}$ of $\R^{n}$ is called an \textit{algebraic interior} associated with $f$ if it is the closure of a (arcwise) connected component of $\{\x \in \R^{n}: f(\x) > 0\}$. The polynomial $f$ is called a \textit{defining polynomial} for $C_{f}.$   Consequently, if $f(\x)$ has an MDR, $C_{f}=S$ \cite{Helton}. But the converse of this statement need not be true \cite{Helton}.

One way to characterize monic (definite) LMI representable sets is by identifying polynomials which have monic (definite) symmetric or Hermitian determinantal representations \cite{Sturmfelsmin}, \cite{Helton}. A recent literature survey in this area can be found in \cite{Vinnikov}. It is to be noted that amongst all spectrahedra, those defined by a LMI (\ref{lmiset}) which have $A_{0} \succ 0$ are special, as problems related to these spectrahedra can be connected to semidefinite programming.

A multivariate polynomial $f(\x) \in \R[\x]$ is said to be a real zero (RZ) polynomial if the polynomial has only real zeros when restricted to any line passing through origin i.e., for any $\x \in \R^{n}$, all the roots of the univariate polynomial $f_{\x}(t):=f(t\cdot \x)$ are real (and $f(0) \neq 0$). If a polynomial $f(\x)$ admits an MDR, say $f(\x)=\det(I+x_{1}A_{1} +x_{2}A_{2} + \dots + x_{n}A_{n})$ then it is indeed a RZ polynomial. This follows from the fact that the univariate polynomial $f_{\x}(t):=\det(I+t(x_{1}A_{1} +x_{2}A_{2} + \dots + x_{n}A_{n}))$ has only real zeros which are in fact the negatives of the reciprocals of non-zero eigenvalues of the Hermitian or symmetric matrix $x_{1}A_{1} +x_{2}A_{2} + \dots + x_{n}A_{n}$ for any given $\x \in \R^n.$ It has been proved that any RZ bivariate polynomial always has an MDR \cite{Helton}. However, if the number of variables of a RZ polynomial is more than $2$, it may not have an MDR at all. For example, dehomogenized polynomial of \textit{Vamos cube} $V_{8}$ is a RZ polynomial without a definite determinantal representation \cite{Branden}.
\section{Polynomials with Monic Determinantal Representations of size $2$} \label{secmsdr2}
In this section, our aim is to characterize all sets which are $2 \times 2$ monic LMI representable and to identify all quadratic polynomials which admit an MDR of size $2$. As RZ property is a necessary condition for a polynomial to have an MDR, we begin with RZ quadratic polynomials.
\subsection{RZ property for Quadratic Polynomials}
It is well known that any quadratic polynomial $f(\x) \in \R[\x]$ (where $\x=(x_{1}, \dots, x_{n})$) can be written
as $f = Z^{T}[\x] Q Z[\x]$ where $Q \in \R^{(n+1) \times (n+1)}$
 and $Z[\x] =\bmatrix{1 & x_{1} & \dots & x_{n}}^{T}$.
Such a matrix $Q$ associated with $f$ is unique if $Q \in S\R^{(n+1) \times (n+1)}$, -- this matrix $Q$ is referred to as the matrix representation of the polynomial $f(\x)$.
The following proposition provides a necessary and sufficient condition for a quadratic polynomial to be a RZ polynomial that shall be used in sequel.
\begin{proposition} \label{rznasc}
Let $f(\x)=Z^{T}[\x] Q Z[\x] \in \R[\x]$ be a quadratic polynomial with nonzero constant term $c$, and let $Q \in S\R^{(n+1) \times (n+1)}$ be the matrix representation of $f(\x)$. Then $f(\x)$ is a RZ polynomial if and only if the Schur complement of $Q$ with respect to the $(1,1)$ element $c$ of $Q$ is negative semidefinite.
\end{proposition}
\pf Note that $f(\x) \in \R[\x]$ can be written as $f(\x) = \x^{T}A\x + b^{T}\x +c,$ where $c=f(0) \neq 0$. Then 
\begin{equation} \label{matrixrepresentation}
f(\x)= \left[\begin{array} {cc}
     1 & \x
       \end{array}\right] \left[\begin{array} {cc}
     c & b^{T}/2 \\
     b/2 & A
       \end{array}\right] \left[\begin{array} {c}
     1 \\
     \x \\
       \end{array}\right] = Z^{T}[\x]QZ[\x].
\end{equation}
The Schur complement of $(1,1)$ element $c$ in the matrix representation $Q$ is $A- \frac{b}{2}(1/c)\frac{b^{T}}{2}$.
As $f(\x) = \x^{T}A\x + b^{T}\x +c$, so $f_{\x}(t)=f(t\x) = t^{2}\x^{T}A\x + t b^{T}\x +c$, $t\in \R$. Therefore, for any $\x \in \R^{n},$
\begin{align*}
&\mbox{the roots of polynomial} \ f_{\x}(t) \ \mbox{ are real}    \\
\Leftrightarrow &(b^{T}\x)^{2} - 4c \x^{T}A\x \geq 0 \\
\Leftrightarrow &\x^{T}(bb^{T}-  4Ac)\x \geq 0  \\
&\Leftrightarrow 4Ac - bb^{T} \preceq 0.
 \end{align*} \qed

Any polynomial admitting an MDR is a RZ polynomial, but the converse need not be true. We give an example below where the converse does not hold.
\begin{example} \cite{Tim} \label{nrzquad}
The (shifted hyperbolic) polynomial $f(\x)=(x_{1}+1)^{2} - x_{2}^{2} - x_{3}^{2} - x_{4}^{2}$ is a quadratic RZ polynomial which has no MSDR, but it has an MHDR of size $2$. On the other hand the polynomial $f(\x)=1-x_{1}^{2}-x_{2}^{2}-x_{3}^{2}-x_{4}^{2}$ has no MHDR (so it can not have an MSDR too).
\end{example}
We now provide a necessary and sufficient condition for the existence of MDR of size $2$ for RZ quadratic polynomials.

\subsection{Quadratic polynomials having MDR of size $2$}
In this subsection, we provide a necessary and sufficient condition for a quadratic polynomial to have an MDR of size $2$. We further provide an algorithm to determine MDRs for a quadratic polynomial when they exist. Note that since we are interested in monic representations, therefore the constant term of the quadratic polynomial must be one. Henceforth, for a quadratic polynomial $f(\x) \in \R[\x]$ we denote the Schur complement of its matrix representation $Q$ with respect to the $(1,1)$ element by $Q/(1,1)$.
\begin{theorem} \label{sufficientqmsdr}
A quadratic polynomial $f(\x)=Z^{T}[\x]QZ[\x] \in \R[\x]$ with $f(0)=1$ has an MDR of size $2$
if and only if $Q/(1,1)$ is negative semidefinite and the rank$(Q/(1,1)) \leq 3$. If $Q/(1, 1)$ is negative semidefinite and the rank$(Q/(1,1)) = 3$, then the polynomial has Hermitian MDR but no symmetric MDR. For symmetric MDR, $Q/(1,1)$ must be negative
semidefinite and rank$(Q/(1, 1)) \leq 2$.
\end{theorem}
\pf Suppose $f(\x)=\det(I+x_{1}A_{1} +x_{2}A_{2} + \dots + x_{n}A_{n})$ where
\begin{equation} \nonumber
A_{j}= \left[\begin{array} {cc}
     r_{j} & t_{j}-iu_{j}  \\
     t_{j}+iu_{j} & s_{j}
       \end{array}\right] \in \H^{2 \times 2}(\C), j=1,\dots,n.
\end{equation}
Let
\begin{equation*}
\rq=\bmatrix{1 & r_{1} & \dots & r_{n}}^{T}, \s=\bmatrix{1 & s_{1} & \dots & s_{n}}^{T}, \t=\bmatrix{0 & t_{1} & \dots & t_{n}}^{T}, \u=\bmatrix{0 & u_{1} & \dots & u_{n}}^{T}
\end{equation*}
Consider the truncated vectors
\begin{equation*}
\tilde{\rq}=\bmatrix{r_{1} & \dots & r_{n}}^{T}, \tilde{\s}=\bmatrix{s_{1} & \dots & s_{n}}^{T}, \tilde{\t}=\bmatrix{t_{1} & \dots & t_{n}}^{T}, \tilde{\u}=\bmatrix{u_{1} & \dots & u_{n}}^{T}
\end{equation*}
of $\rq,\s,\t$ and $\u$ respectively.
Note that 
\begin{equation} \nonumber
Q = (1/2)(\rq \s^{T} + \s \rq^{T}) + (-\t \t^{T})+ (-\u \u^{T}).
\end{equation}
Therefore, rank$(Q)\leq 4$, as $Q$ is the sum of four rank one matrices. Consequently, $Q/(1,1)$ is given by
\begin{align*}
&(1/2)(\tilde{\rq}\tilde{\s}^{T} + \tilde{\s}\tilde{\rq}^{T}) + (-\tilde{\t} \tilde{\t}^{T})+(-\tilde{\u} \tilde{\u}^{T})+ [- (\frac{\tilde{\rq}+\tilde{\s}}{2}) (\frac{\tilde{\rq}+\tilde{\s}}{2})^{T}]\\
&= -((\frac{\tilde{\rq}-\tilde{\s}}{2}) (\frac{\tilde{\rq}-\tilde{\s}}{2})^{T})+ (-\tilde{\t} \tilde{\t}^{T})+ (-\tilde{\u} \tilde{\u}^{T}).
\end{align*}
Thus rank$(Q/(1,1)) \leq 3$ follows from the fact that $Q/(1,1)$ is the sum of three rank one matrices. It is known that if a quadratic polynomial $f(\x) \in \R[\x]$ has an MDR, then $f(\x)$ is a RZ polynomial. By Proposition \ref{rznasc}, $Q/(1,1)$  is negative semidefinite. Hence, if a quadratic polynomial $f(\x)$ with $f(0)=1$ has an MDR of size $2$, then rank$(Q/(1,1))\leq 3$ and $Q/(1,1)$ is a negative semidefinite matrix.

Conversely, suppose $f(\x) \in \R[\x]$ is a quadratic polynomial with $f(0)=1$ and $Q/(1,1)$ is a negative semidefinite matrix such that rank$(Q/(1,1)) \leq 3$.
Since $Q/(1,1)$ is the Schur complement with respect to the $(1, 1)$ element of the matrix
$Q=\bmatrix{1 & b^{T}/2 \\ b/2 & A}$ which represents the quadratic polynomial $f(\x)=\x^{T}A \x+b^{T} \x +1$, 
so we have $Q/(1,1)= A-(1/4)bb^{T}$.
Thus to obtain an MDR we need to find the vectors $\tilde{\rq},\tilde{\s},\tilde{\t}$ and $\tilde{\u}$, which were defined earlier.
Since rank of $(Q/(1,1))$ is at most $3$, we can obtain $-(Q/(1,1))$ as sum of three rank one matrices $\alpha_{1}\alpha_{1}^{T}+\alpha_{2}\alpha_{2}^{T}+\alpha_{3}\alpha_{3}^{T}$ (for example, by using the \textit{Cholesky decomposition}).
In fact,
\begin{equation} \nonumber
-(Q/(1,1))=(1/4)b b^{T}-A=((\frac{\tilde{\rq}-\tilde{\s}}{2}) (\frac{\tilde{\rq}-\tilde{\s}}{2})^{T}) +(\tilde{\t} \tilde{\t}^{T})+(\tilde{\u} \tilde{\u}^{T})
\end{equation}
Note that $\tilde{\rq}+\tilde{\s}=b$ where $b$ is defined in the matrix representation $Q$. Therefore one can obtain a Hermitian MDR by setting $(1/2)(\tilde{\rq}-\tilde{\s})=\alpha_{1}, \tilde{\t}=\alpha_{2}$, and $\tilde{\u}=\alpha_{3}$, and solving this system of linear equations along with $\tilde{\rq}+\tilde{\s}=b$.

Note that for a symmetric MDR, the coefficient matrices $A_{j}$ are of the form $\bmatrix{r_{j} & t_{j}\\t_{j} & s_{j}}$ which make $A_{j}$s symmetric. Hence from the proof above, it is clear that the vector $\tilde{\u}$ must now be the zero vector.
Thus $Q/(1,1) $ must be the sum of two rank one matrices and therefore the rank of matrix $Q/(1,1) \leq 2$. \hfill{$\square$}

As a consequence of the above theorem we have the following corollary which provides a necessary and sufficient condition for a quadratic polynomial to have a diagonal MDR of size $2$.
\begin{corollary} \label{diagmsdrcor}
A quadratic polynomial $f(\x)=Z^{T}[\x]QZ[\x], f(0)=1$ has a diagonal MDR of size $2$ if and only if $Q/(1,1)$ is negative semidefinite and of rank at most $1$.
\end{corollary}
\pf The proof follows from the proof of the above theorem by setting $\tilde{\t}=0$ and $\tilde{\u}=0$. \hfill${\square}$

Based on the constructive proof of Theorem \ref{sufficientqmsdr}, we now provide an algorithm to determine MDR of size $2$ for a quadratic polynomial, whenever it exists.
\begin{algorithm} \label{algo}
\caption{Algorithm to find MHDR (MSDR) of size $2$}
\begin{algorithmic}
\State Input: a quadratic polynomial $f(\x)=\x^{T}A\x+b^{T}\x+1$ 
\State Output: $A_{j}=\bmatrix{r_{j} & t_{j}-iu_{j} \\ t_{i}+iu_{j} & s_{i}}, j=1,\dots ,n$ such that $f(\x)=\det\left(I+x_{1}A_{1}+\dots+x_{n}A_{n}\right)$.
\\\hrulefill
\begin{enumerate}
\item Form the matrix representation $Q$ of $f(\x)$
\item Calculate $Q/(1,1)=A-\frac{1}{4} bb^{T}$
\item Compute the Cholesky factor of $-Q/(1,1)$. If Cholesky factor does not exist, no MHDR (MSDR) of size $2$ possible.
\item If rank$(Q/(1,1)) > 3$, then exit -- no MDR of size $2$ possible 
\item Otherwise, $Q/(1,1)=\alpha_{1}\alpha_{1}^{T}+\alpha_{2}\alpha_{2}^{T}+\alpha_{3}\alpha_{3}^{T}$. 
\item If rank $(Q/(1,1)) = 3$, then construct a Hermitian MDR of size $2$ by setting $r_{j}=(\frac{(2\alpha_{1}+b)}{2})_{j}, s_{j}=(\frac{(b-2\alpha_{1})}{2})_{j}, t_{j}=(\alpha_{2})_{j}, u_{j}=(\alpha_{3})_{j}$.
\item Construct $A_{j}=\bmatrix{r_{j} & t_{j}-iu_{j} \\ t_{j}+iu_{j} & s_{j}}$ for $j=1, \dots,n$.
\item If rank $(Q/(1,1)) \leq 2$, then obtain a symmetric MDR by setting $u_{j}=0, r_{j}=(\frac{(2\alpha_{2}+b)}{2})_{j}, s_{j}=(\frac{(b-2\alpha_{2})}{2})_{j}$, and $t_{j}=(\alpha_{1})_{j}$ 
\item Assign $A_{j}=\bmatrix{r_{j} & t_{j} \\ t_{j} & s_{j}}$ for $j=1, \dots,n$.
\end{enumerate}
\end{algorithmic}
\end{algorithm}

We demonstrate the algorithm in the following examples.
 \begin{example}\label{exquniequiv} {\rm  We provide three examples below.
\begin{enumerate}
\item Consider $f(\x)= 1-8 x_{1} x_{2} -4 x_{1}x_{3} -100 x_{2}^{2}-12 x_{2} x_{3} -x_{3}^{2} -5 x_{1}^{2}$. Then $f(\x)=Z^T[\x]Q Z[\x]$ where $$ Q =\bmatrix{1 & 0 & 0 & 0 \\ 0 & -5 & -4 & -2 \\0 & -4 & -100 & -6 \\ 0 & -2 & -6 & -1}.$$ Then it is easy to verify that $Q/(1,1)$ is negative semidefinite and the rank of $Q/(1,1)$ is $2.$ Indeed,
   \begin{align*}
     -Q/(1,1) &=\bmatrix{5 & 4 & 2 \\4 & 100 & 6 \\2 & 6 & 1}
     =\bmatrix{2/5 \\ 10 \\ 3/5} \bmatrix{2/5 & 10 & 3/5} + \bmatrix{11/5 \\ 0 \\ 4/5}\bmatrix{11/5 & 0 & 4/5}
    \end{align*}
    Hence $f(\x)$ admits an MSDR. By applying Algorithm $1$, one obtains
 \begin{align*}
f(\x) =\det \left( I+x_{1}\bmatrix{11/5 & 2/5 \\ 2/5 & -11/5} + x_{2}\bmatrix{0 & 10 \\ 10 & 0} + x_{3} \bmatrix{4/5 & 3/5 \\ 3/5 & -4/5} \right),
\end{align*}
\item Consider $f(\x)= 1+4x_{1}-10x_{2}-x_{1}^{2}-2x_{1}x_{2}-x_{2}^{2}$ which has a Hermitian determinantal representation provided in the chapter \cite{Vinnikov2}. In fact $f(\x)$ admits a symmetric MDR, since $Q=\bmatrix{1 & 2 & 5 \\ 2 & -1 & -1 \\5 & -1 & -1}$ and $Q/(1,1)$ is negative semidefinite of rank $2$.
     \begin{equation*}
-Q/(1,1)=\bmatrix{5 & 11 \\ 11 & 26}=\bmatrix{2.2316 \\
    4.9193}\bmatrix{2.2316 & 4.9193} + \bmatrix{0 \\ 1.3416}\bmatrix{0 & 1.3416}.
    \end{equation*}
Then by Algorithm $1$,
\begin{align*}
f(\x) =\det\left(I+x_{1}\bmatrix{2 & 2.2361 \\ 2.2361 & 2}+x_{2}\bmatrix{6.3416 & 4.9193 \\4.9193 & 3.6584}\right).
\end{align*}
\item Consider the quadratic polynomial $f(\x)=(x_{1}+1)^{2}-x_{2}^{2}-x_{3}^{2}-x_{4}^{2}$. Here matrix representation
\begin{equation*}
Q=\bmatrix{1 & 1 & 0 & 0 &0\\1 & 1 & 0 & 0 & 0\\0 & 0 & -1 & 0 & 0\\0 & 0 & 0 & -1 & 0\\0 & 0 & 0 & 0 & -1}, -Q/(1,1)=\bmatrix{0 & 0 & 0 & 0\\0 & 1 & 0 &0\\0& 0 & 1 & 0\\0 & 0 & 0 &1}=\alpha_{1}\alpha_{1}^{T}+\alpha_{2}\alpha_{2}^{T}+\alpha_{3}\alpha_{3}^{T}.
\end{equation*}
where one could take $\alpha_{1}=\bmatrix{0\\0\\0\\1},\alpha_{2}=\bmatrix{0\\0\\1\\0},\alpha_{3}=\bmatrix{0\\1\\0\\0}.$
Since the rank of the matrix is $3$, there is no symmetric MDR, but then Hermitian MDR does exist. Coefficient
matrices for a Hermitian MDR are as follows.
$$A_{1}=\bmatrix{1 & 0\\0 & 1}, A_{2}=\bmatrix{0 & -i\\i & 0}, A_{3}= \bmatrix{0 & 1\\1 & 0}, A_{4}=\bmatrix{1 & 0\\0 & -1}.$$
\end{enumerate}
}
\end{example}
\subsection{Equivalent and Non-equivalent MSDRs and MHDRs}
Note that in the algorithm given above, there is considerable freedom in constructing the
MDR. For one, the choice of the vectors $\alpha_{1}, \alpha_{2},\alpha_{3}$ to express $Q/(1,1)$ as the sum of three
rank one matrices is immense. Each such choice, leads to a different MDR.
Two linear matrix polynomials $A_{0}+x_{1}A_{1}+\dots+x_{n}A_{n}$ and $B_{0}+x_{1}B_{1}+\dots+x_{n}B_{n}$ of size
$k$ are said to be unitarily (orthogonally) equivalent if there exists an unitary (orthogonal)
matrix $U$ of order $k$ such that $U(A_{0} + x_{1}A_{1} + \dots + x_{n}A_{n})U^{\ast} = B_{0} + x_{1}B_{1} + \dots + x_{n}B_{n}$.
Note that for an orthogonal matrix $U^{\ast} = U^{T}$ . Note that the determinants of two equivalent
linear matrix polynomials are the same. Using this equivalence, we can declare two different MDRs of a given polynomial as equivalent, if the corresponding linear matrix polynomials are equivalent. The equivalence class of representations is said to be definite (monic) if it
contains a definite (monic) representative. Naturally, one would be interested in determining how many non-equivalent classes of MDR exists for a given polynomial. Observe that the matrices $A_{j}$ obtained in the algorithm given above, can be re-written as the sum of a diagonal matrix and a traceless matrix. Thus
\begin{equation}
A_{j}=\frac{b_{j}}{2}\bmatrix{1 & 0 \\ 0 & 1}+\bmatrix{(\alpha_{1})_{j} & (\alpha_{2})_{j}-i(\alpha_{3})_{j} \\ (\alpha_{2})_{j}+i(\alpha_{3})_{j} & -(\alpha_{1})_{j}}
\end{equation}
Here $b_{j}, (\alpha_{1})_{j}, (\alpha_{2})_{j},(\alpha_{3})_{j}$ are the $j$-th entries of the vectors $b, \alpha_{1}, \alpha_{2},\alpha_{3}$ respectively. Similarly, for the symmetric case, $A_{j}=\frac{b_{j}}{2}\bmatrix{1 & 0 \\ 0 & 1}+\bmatrix{(\alpha_{1})_{j} & (\alpha_{2})_{j} \\ (\alpha_{2})_{j} & -(\alpha_{1})_{j}}$. Note further that similarity transforms using unitary (orthogonal) matrices does not affect the diagonal matrix, i.e.,
$U^{\ast}A_{j}U=\frac{b_{j}}{2}\bmatrix{1 & 0 \\ 0 & 1}+U^{\ast}\bmatrix{(\alpha_{1})_{j} & (\alpha_{2})_{j}-i(\alpha_{3})_{j} \\ (\alpha_{2})_{j}+i(\alpha_{3})_{j} & -(\alpha_{1})_{j}}U$. Thus, it is enough to consider unitary (orthogonal) equivalence of traceless Hermitian (symmetric) matrices. Traceless Hermitian matrices form a three dimesional real vector space spanned by the Pauli matrices $\sigma_{z}=\bmatrix{1 & 0\\0 & -1}, \sigma_{x}=\bmatrix{0 & 1\\1 & 0}, \sigma_{y}=\bmatrix{0 & -i\\i & 0}$. Similarly, traceless symmetric matrices form a
two dimensional real vector space spanned by $\sigma_{z}=\bmatrix{1 & 0\\0 & -1}$ and $\sigma_{x}=\bmatrix{0 & 1\\1 & 0}$.  Therefore, to
figure out equivalent Hermitian MDRs, one needs to consider only equivalence of Hermitian matrices of the form $(\alpha_{1})_{j}\sigma_{z}+(\alpha_{2})_{j}\sigma_{x}+(\alpha_{3})_{j}\sigma_{y}$. Similarly, for symmetric MDRs, it is enough
to consider equivalence of symmetric matrices of the form $(\alpha_{1})_{j}\sigma_{z}+(\alpha_{2})_{j}\sigma_{x}$. We first consider the case of symmetric MDR. Consider an orthogonal matrix $V=\bmatrix{v_{11} & v_{12}\\v_{21} & v_{22}}$. One can compute that
\begin{equation} \label{eqqequiv1}
V^{T}(k \sigma_{z}+l \sigma_{x})V=((v_{11}^{2}-v_{21}^{2})k+2v_{11}v_{21}l)\sigma_{z}+((v_{11}v_{12}-v_{21}v_{22})k+(v_{11}v_{22}+v_{12}v_{21})l)\sigma_{x}
\end{equation}
We have shown that a symmetric MDR of size $2$ for a quadratic polynomial exists if and
only if the Schur complement matrix $Q/(1, 1)$ is negative semidefinite with rank less than
or equal to $2$. Let $-Q/(1,1) = R^{T}R = \alpha_{1}\alpha_{1}^{T}+\alpha_{2}\alpha_{2}^{T}$ which can be used to construct the
coefficient matrices $A_{j}$s as outlined in the earlier section. Here $\alpha_{1}^{T},\alpha_{2}^{T}$ are the rows of the
matrix $R \in \R^{2 \times n}$. Note that we can construct two alternate symmetric MDRs by either
assigning $\alpha_{1}=\tilde{\t}, \alpha_{2}=\frac{\tilde{\rq}-\tilde{\s}}{2}$ or $\alpha_{1}=\frac{\tilde{\rq}-\tilde{\s}}{2}, \alpha_{2}=\tilde{\t}$.
\begin{proposition}
All symmetric MDRs of size $2$ of a quadratic polynomial are orthogonally equivalent.
\end{proposition}
\pf Let $-Q/(1, 1) = R^{T}R$ where $R \in \R^{2 \times n}$. If the $j$-th column of $R$ is given by
the vector $\bmatrix{k \\l}$, then the matrix $A_{j}$ of the corresponding MDR is by construction equal
to $\frac{b_{j}}{2} I + k\sigma_{z} + l \sigma_{x}$. If $-Q/(1,1)= (OR)^{T}OR$ where $O$ is any orthogonal $2 \times 2$ matrix,
then the corresponding matrix of the new MDR is equal to $\frac{b_{j}}{2} I + k_{1}\sigma_{z} + l_{1} \sigma_{x}$ where $\bmatrix{k_{1}\\ l_{1}}=O\bmatrix{k\\l}$. Observe that if $O =\bmatrix{v_{11}^{2}-v_{21}^{2} & 2v_{11}v_{21}\\v_{11}v_{12}-v_{21}v_{22} & v_{11}v_{22}+v_{12}v_{21}}$, then by equation (\ref{eqqequiv1}) above, the matrices associated to the new MDR are equal to $V^{T}A_{j}V$. Thus, it is enough to demonstrate
that for every orthogonal matrix $O$, there exists another orthogonal matrix $V$ such that the elements of $O$ are related to those of $V$ in the manner described above. If the orthogonal matrix $O$ is a rotation matrix: \ $O=\bmatrix{\cos \theta & \sin \theta\\-\sin \theta & \cos \theta}$, then $V=\bmatrix{\cos \frac{\theta}{2} & -\sin \frac{\theta}{2}\\\sin \frac{\theta}{2} & \cos \frac{\theta}{2}}$ satisfies the required relations. On the other hand, if the orthogonal matrix $O$ is a reflection matrix $O=\bmatrix{-\cos \theta & -\sin \theta\\-\sin \theta & \cos \theta}$, then $V=\bmatrix{-\sin \frac{\theta}{2} & \cos \frac{\theta}{2}\\ \cos \frac{\theta}{2} & \sin \frac{\theta}{2}}$ satisfies the required relations. \hfill{$\square$}

We now consider the case of Hermitian MDRs. We begin with some remarks about
unitary matrices. The set of unitary matrices $U(2)$ is a Lie group and consists of matrices $U$ such that $U^{\ast}U = I_{2}$. Clearly $\det(U) = e^{i \phi}$, for $U \in U(2)$. The set of $2 \times 2$ unitary matrices with determinant $1$ is also a Lie group, denoted by $SU(2)$. Given $U \in U(2)$ with determinant
$e^{i\phi}$, observe that the matrix $U_{1} = e^{\frac{-i \phi}{2}}U \in SU(2)$. Further, for any Hermitian matrix $A$, $U^{\ast}AU = U_{1}^{\ast}AU_{1}$. Thus, for unitary equivalence of Hermitian MDRs, it is enough to consider unitary matrices from $SU(2)$.

Any $U \in SU(2)$ can be written as $U = \bmatrix{a+ib &  -c + id \\ c+ id & a- ib}$ where $a^{2} + b^{2} + c^{2} + d^{2} = 1$.
If $H$ is a traceless Hermitian matrix, then $H_{1} = U^{\ast}HU$ is also a traceless Hermitian matrix.
Further, if one expresses these traceless Hermitian matrices in terms of the Pauli matrices as
$H = k\sigma_{z} + l\sigma_{x} + m\sigma_{y}$ and $H_{1} = k_{1}\sigma_{z} + l_{1}\sigma_{x} + m_{1}\sigma_{y}$, then
\begin{equation} \label{eqqequiv2}
\bmatrix{k_{1}\\l_{1}\\m_{1}}=\bmatrix{a^{2}+b^{2}-c^{2}-d^{2} & 2ac+2bd & 2ad-2bc \\2bd-2ac & a^{2}-b^{2}-c^{2}+d^{2} & -2ab-2cd\\-2ad-2bc & 2ab-2cd & a^{2}-b^{2}+c^{2}-d^{2}} \bmatrix{k\\l\\m}
\end{equation}

We have shown that a Hermitian MDR of size $2$ exists for a quadratic polynomial if and
only if $Q/(1, 1)$ is negative semidefinite with rank $(Q/(1, 1)) \leq 3$. Let $-Q/(1,1) = R^{T}R =\alpha_{1}\alpha_{1}^{T}+\alpha_{2}\alpha_{2}^{T}+\alpha_{3}\alpha_{3}^{T}$ where $\alpha_{1}^{T},\alpha_{2}^{T},\alpha_{3}^{T}$ are the rows of $R \in \R^{3 \times n}$.

\begin{proposition}
Quadratic polynomials that have a Hermitian MDR of size $2$ but no symmetric MDR have two classes of unitarily equivalent MDRs. All other quadratic polynomials that have a MDR of size $2$ have only one class of unitarily equivalent MDRs.
\end{proposition}
\pf By Theorem \ref{sufficientqmsdr}, if a quadratic polynomial has a Hermitian MDR of size $2$ but no symmetric MDR, then $-Q/(1, 1) = R^{T}R$ where the full row rank matrix $R \in \R^{3 \times n}$. Observe that if the $j$-th column of $R$ is given by $\bmatrix{k\\l\\m}$, then the matrix $A_{j}$ of the MDR
is given by $\frac{b_{j}}{2} I+k \sigma_{z}+l \sigma_{x}+m \sigma_{y}$. One can obtain another factorization of $-Q/(1, 1)$ as
$-Q/(1, 1) = (OR)^{T}OR$ where $O$ is a $3 \times 3$ orthogonal matrix. The MDR obtained from this new factorization would have $A_{j}=\frac{b_{j}}{2} I+k_{1} \sigma_{z}+l_{1} \sigma_{x}+m_{1} \sigma_{y}$ where $\bmatrix{k_{1}\\l_{1}\\m_{1}}=O \bmatrix{k\\l\\m}$. If $O$
has the form of the matrix in equation \ref{eqqequiv2}, then the MDRs obtained by the two factorizations
are unitarily equivalent. The determinant of the matrix in equation \ref{eqqequiv2} is $1$ and therefore
determinant of $O$ must be $1$ for the two MDRs to be unitarily equivalent. As an orthogonal
matrix can have determinant equal to $\pm 1$, therefore if one uses an orthogonal matrix $O$ with
determinant equal to $-1$, then the two MDRs obtained from $R$ and $OR$ are not unitarily
equivalent. Thus there are two classes of unitarily equivalent MDRs.

Now we consider the case of a quadratic polynomial where rank$(Q/(1, 1)) < 3$. In this
case, $-Q/(1, 1) = R^{T}R$ where the full row rank matrix $R \in \R^{s \times n}$ with $s < 3$. We can also
view this as $-Q/(1, 1) = R^{T}_{1} R_{1}$ where $R_{1} \in \R^{3 \times n}$ has been obtained from $R$ by appending all zero row(s). In this particular case, one can obtain equivalent Hermitian MDRs by modifying $R_{1}$ to $OR_{1}$, where $O$ is a $3 \times 3$ orthogonal matrix where the all-zero rows are preserved as all-zero rows of the new matrix. Such a transformation ensures that $-Q/(1, 1) = R^{T}_{1} R_{1} =(OR_{1})^{T}OR_{1}$.

We therefore explore what happens when one equates an $O$ that preserves the all-zero
rows of $R_{1}$ to the matrix from equation \ref{eqqequiv2}
\begin{equation*}
\bmatrix{a^{2}+b^{2}-c^{2}-d^{2} & 2ac+2bd & 2ad-2bc \\2bd-2ac & a^{2}-b^{2}-c^{2}+d^{2} & -2ab-2cd\\-2ad-2bc & 2ab-2cd & a^{2}-b^{2}+c^{2}-d^{2}}
\end{equation*}
Let us assume that the rank$(-Q/(1, 1)) = 2$ and $R_{1}$ is a matrix whose third row is the all
zero row. Therefore one of the rows of the orthogonal matrix $O$ must be $\pm e_{3} = (0, 0,\pm 1)$ to preserve the all zero row. If the third row is $e_{3}$, then $a^{2} + c^{2} = 1$ and $b = d = 0$. This automatically ensures that all the other elements of the third row and third column are zero.
Thus we obtain $O$ to have the form
\begin{equation*}
\bmatrix{a^{2}-c^{2} & 2ac & 0\\ -2ac & a^{2}-c^{2} & 0\\0 & 0 & a^{2}+c^{2}}
\end{equation*}
Thus we observe that a rotation matrix is applied to the first two rows of the matrix $R_{1}$. On
the other hand, if the third row is $\pm e_{3}$, then $b^{2} + d^{2} = 1$ and $a = c = 0$ giving the matrix
\begin{equation*}
\bmatrix{b^{2}-d^{2} & 2bd & 0\\ 2bd & -b^{2}+d^{2} & 0\\ 0 & 0 & -b^{2}-d^{2}}
\end{equation*}
This corresponds to the reflection matrix being applied to the first two rows of $R_{1}$.

Now consider the case where the first row of the matrix $O$ is $\pm e_{3}$. For such a $O$ to be of the
form given by equation (\ref{eqqequiv2}, equating the expressions from the first row and the third column,
one obtains either $a = d, b = -c$ with $a^{2} + b^{2} = 1/2$ when the first row is $e_{3}$ or $a = -d,b = c$
with $a^{2} + b^{2} = 1/2$ when the first row is $-e_{3}$. For the first case, the $2 \times 2$ submatrix of $O$
acting on the nontrivial rows of $R_{1}$ is a rotation matrix whereas for the second case, this
submatrix is a reflection matrix. Similarly, if one assumes the second row of $O$ is $\pm e_{3}$, then
one can show that the relevant $2 \times 2$ submatrix of $O$ that acts on the nontrivial rows of $R_{1}$ is
a reflection matrix when the second row of $O$ is $e_{3}$ whereas it is a rotation matrix when the
second row of $O$ is $-e_{3}$.

This clearly shows that all unitarily equivalent MDRs for this case are indeed obtained
from the original $R_{1}$ by a $3 \times 3$ orthogonal matrix having determinant equal to $1$. The case
where $R_{1}$ has only one nontrivial row is trivial, since that row or its negative should be the
only nontrivial row of $OR_{1}$ and this is easily obtained with $O$ having determinant equal to
$1$. \hfill{$\square$}

We can demonstrate this by using the earlier Example \ref{exquniequiv}.
\begin{example} \label{exquniequiv2} \rm{
\begin{enumerate}
\item Recall $1$ of Example \ref{exquniequiv}.
The sets of coefficient matrices
\begin{equation*}
\{\bmatrix{11/5 & 2/5 \\ 2/5 & -11/5},\bmatrix{0 & 10 \\ 10 & 0}, \bmatrix{4/5 & 3/5 \\ 3/5 & -4/5}\} \mbox{and }\{\bmatrix{11/5 & 2i/5 \\ -2i/5 & -11/5},\bmatrix{0 & 10i \\ -10i & 0}, \bmatrix{4/5 & 3i/5 \\ -3i/5 & -4/5}\}
 \end{equation*}
are unitarily equivalent by matrix $U=\bmatrix{\frac{1-i}{\sqrt{2}} & 0 \\ 0 & \frac{1+i}{\sqrt{2}}}$.
Another unitarily equivalent
MDR is given by the coefficient matrices
\begin{equation*}
\bmatrix{4/5 & 3/5 \\ 3/5 & -4/5}\} \mbox{and} \{\bmatrix{2/5 & 11i/5 \\ -11i/5 & -2/5},\bmatrix{10 & 0 \\ 0 & -10}, \bmatrix{3/5 & 4i/5 \\ -4i/5 & -3/5}\}
\end{equation*}
and these are obtained by using unitary matrix $U=\bmatrix{\frac{1+i}{2} & -\frac{1-i}{2} \\ \frac{1+i}{2} & \frac{1-i}{2}}$ on the original coefficient matrices.
\item From part $3$ of Example \ref{exquniequiv}, recall that coefficient matrices for a Hermitian MDR for
the polynomial $f(\x)=(x_{1}+1)^{2}-x_{2}^{2}-x_{3}^{2}-x_{4}^{2}$ were
\begin{equation*}
A_{1}=\bmatrix{1 & 0\\0 & 1}, A_{2}=\bmatrix{0 & -i\\i & 0}, A_{3}= \bmatrix{0 & 1\\1 & 0}, A_{4}=\bmatrix{1 & 0\\0 & -1}.
\end{equation*}

Observe that another Hermitian MDR for the same polynomial is given by the coefficient
matrices
\begin{equation*}
A_{1}=\bmatrix{1 & 0\\0 & 1}, A_{2}=\bmatrix{0 & 1\\1 & 0}, A_{3}= \bmatrix{0 & -i\\i & 0}, A_{4}=\bmatrix{1 & 0\\0 & -1}.
\end{equation*}
These two MDRs are not unitarily equivalent. This is because the vectors $\alpha_{2}$ and $\alpha_{3}$
were swapped in the original factorization to obtain the second MDR from the first
one. This swapping of vectors arises out of the action of a permutation matrix whose
determinant is $-1$.
\end{enumerate}
}
\end{example}
\section{Complete characterization of quadratic polynomials that admit MDRs} \label{seccompmsdr}
We now completely characterize all quadratic polynomials which are determinants of monic linear matrix polynomials of any size.
\subsection{Spectrahedra}
We assume that a spectrahedron has a nonempty interior and therefore without loss of generality, we assume that the spectrahedron contains origin as an interior point. Therefore it is determined by a definite LMP \cite{Ramana1}, \cite{Tim}. Now we define what is meant by the expression ``a spectrahedron contains a full dimensional cone'' \cite{Tim}.
\begin{definition}
Consider a spectrahedron $S =\{\x \in \R^{n}: I+x_{1}A_{1}+\dots+x_{n}A_{n} \succeq 0\}$. Let $f(\x) =\det(I+x_{1}A_{1}+\dots+x_{n}A_{n})$ and $d$ is the degree of the polynomial $f(\x)$. Then the spectrahedron $S$ contains a full dimensional cone if and only if the half ray through some point $\x \in \R^{n}$ is contained within the spectrahedron $S$ and rank$(x_{1}A_{1} +x_{2}A_{2} + \dots + x_{n}A_{n}) = d$. 
\end{definition}

Given a point $\x \in \R^n$, the half ray through the point $\x$ is the set of points obtained as $\{ \lambda \x : \lambda \geq 0 \}$.
 Observe that $x_{1}A_{1}+\dots+x_{n}A_{n} \succeq 0$ if and only if $\lambda x_{1}A_{1}+\dots+ \lambda x_{n}A_{n} \succeq 0$ for every $\lambda \geq 0$. So,
 a spectrahedron $S=\{\x \in \R^{n}: I+x_{1}A_{1}+\dots+x_{n}A_{n} \succeq 0\}$ contains a full dimensional cone if and only if there exists some $\x=(x_{1},\dots,x_{n}) \in \R^{n}$ such that $x_{1}A_{1}+\dots+x_{n}A_{n} \succeq 0$ with rank$(x_{1}A_{1}+\dots+x_{n}A_{n})=d$.
The following theorem \cite{Tim} illustrates the connection between a spectrahedron $S$ containing a full dimensional cone and an MDR of the polynomial $f(\x).$
\begin{theorem}  \label{fulldimensional}
Let $f(\x) =\det(I+x_{1}A_{1}+\dots+x_{n}A_{n}) \in \R[\x]$ where $A_{i} \in S \R^{k \times k}$ or $A_{i} \in \H^{k \times k}(\C)$, for some $k$, and the degree of $f(\x)$ be $d$. If the spectrahedron $S:=\{\x \in \R^{n}: I+x_{1}A_{1} +x_{2}A_{2} + \dots + x_{n}A_{n} \succeq 0\}$ contains a full dimensional cone, then the polynomial $f(\x) \in \R[\x]$ admits an MDR of size $d$.
\end{theorem}
Note that this theorem guarantees the existence of an MDR of size $d$ for a polynomial $f(\x)$ of degree $d$ if there exists an MDR of some size $k$ for the polynomial $f(\x).$ Also note that the converse of the Theorem \ref{fulldimensional} need not be true \cite{Tim}. Notice the polynomial
$f(\x)=1-x_{1}^{2}-x_{2}^{2}$ has a symmetric MDR of size $2$, though the spectrahedron defined by this polynomial does not contain a full dimensional cone.
\subsection{Quadratic Polynomials with MDR of Any Size}
In Section \ref{secmsdr2}, a necessary and sufficient condition for a quadratic polynomial to have an MHDR (MSDR) of size $2$ was provided. Using that result, we now derive a necessary and sufficient condition for the existence of an MDR of any size, for a given quadratic polynomial. Given a quadratic polynomial $f(\x)=\x^{T}A\x+b^{T}\x+1$, where $A$ is negative semidefinite, it is easy to construct a symmetric MDR -- a result well known in literature. On the other hand, the case of $A$ not being negative semidefinite is not well known. We shall throw some light on this case. 
We recall the following proposition.
\begin{proposition} \label{proprzeig}
Let $f:=\det (I + x_{1}A_{1} +x_{2}A_{2} + \dots + x_{n}A_{n}) \in \R[\x]$. Then for each $\x \in \R^{n}$ the nonzero eigenvalues of $x_{1}A_{1} +x_{2}A_{2} + \dots + x_{n}A_{n}$ are in one to one correspondence with the zeros of the univariate polynomial $f_{\x}(t)=f(t \x)$. The correspondence is given by the rule $t \rightarrow  -1 / t$.
\end{proposition}
We prove the main theorem of this section based on the following lemmas that deal with the cases of $A$ not being negative semidefinite.
\begin{lemma} \label{indefinite}
Consider the spectrahedron $S$ defined by a quadratic polynomial $f(\x)=\x^{T}A\x+b^{T}\x+1=Z^{T}[\x]QZ[\x]$, with  $f(\x)=\det(I+x_{1}A_{1} +x_{2}A_{2} + \dots + x_{n}A_{n})$, where $A_{i} \in \H^{k \times k}(\C)$ or $A_{i} \in S\R^{k \times k}$ for some $k \geq 2$. Then $S$ contains a full dimensional cone if $A$ is not a negative semidefinite matrix.
\end{lemma}
\pf If $A$ is not a negative semidefinite matrix, then $A$ has at least one positive eigenvalue, say $\lambda$. Let $v = (v_1,v_2, \cdots, v_n)$ be an eigenvector of $A$ corresponding to that positive eigenvalue $\lambda$. 
Consider the univariate polynomials
 \begin{align*}
 f_{v}(t)=t^{2}v^{T}Av+tb^{T}v+1=t^{2} \lambda ||v||^{2}+tb^{T}v+1 \\
 f_{-v}(t)=t^{2}v^{T}Av-tb^{T}v+1=t^{2} \lambda ||v||^{2}-tb^{T}v+1
\end{align*}
All the coefficients of either $f_{v}(t)$ or $f_{-v}(t)$ are positive. Without loss of generality, let us assume that $b^{T}v >0$ and so all the coefficients of $f_{v}(t)$ are positive. Therefore the polynomial $f_{v}(t)$ has two negative real roots. 
By Proposition~\ref{proprzeig}, $v_{1}A_{1} + v_{2}A_{2} + \dots + v_{n}A_{n}$ is therefore positive semidefinite (its nonzero eigenvalues are the negative reciprocals of the roots of $f_v(t)$). As $f_v(t)$ has only two roots, therefore rank$(v_{1}A_{1} +v_{2}A_{2} + \dots + v_{n}A_{n}) = 2$. Thus the spectrahedron $S$ contains a full dimensional cone.  \hfill${\square}$
\begin{theorem} \label{nsd}
Consider a quadratic polynomial $f(\x)=\x^{T}A\x+b^{T}\x+1=Z^{T}[\x]QZ[\x]$  which an MDR such that $f(\x)=\det(I+L(\x))$, where $L(\x):=x_{1}A_{1}+\dots+x_{n}A_{n}$. Then the spectrahedron $S$ defined by polynomial $f(\x)$ does not contain a full dimensional cone if and only if $A$ is negative semidefinite.
\end{theorem}
\pf The `only if' part of this lemma follows from Lemma \ref{indefinite}.

For proving the `if' part, let us assume that $A$ is negative definite.
So, $\x^{T}A\x < 0$ for any $\x \in \R^{n}\setminus \{0\}$. 
Consider the univariate polynomial
\begin{equation*}
f_{\x}(t)=t^{2}\x^{T}A \x +t b^{T} \x +1.
\end{equation*}
As the coefficient of $t^{2}$ is always negative, the two roots of $f_{\x}(t)$ are real and have opposite signs. 
Therefore using one-to-one correspondence (Proposition \ref{proprzeig}) the non-zero eigenvalues of $L(\x)$ are of opposite signs.
Thus, $L(\x) \nsucceq 0$ for any $\x \in \R^{n}$. Therefore, spectrahedron $S$ does not contain a full dimensional cone.

If $A$ is negative semidefinite, then the way this case differs from the earlier case is that there exists some $u \in \R^{n}\setminus \{0\}$ for which $u^{T}Au=0$. 
This implies the coefficient of $t^{2}$ in $f_{u}(t)$ vanishes, so $f_{u}(t)$ becomes a linear polynomial. Due to Proposition \ref{proprzeig}, the number of non-zero eigenvalues of $L(u)$ is one and so rank$(L(u)) = 1$. Therefore, the spectrahedron does not contain a full dimensional cone, even though it may contain the half ray along $u$. \qed

We demonstrate the construction of a symmetric MDR for a quadratic function $f(\x)$ where $A$ is negative semidefinite. Consider any decomposition (for example, the Cholesky decomposition) of $-A=C^{T}C$. 
This yields
\begin{equation} \label{eqnsd}
f(\x)=-\x^{T}C^{T}C \x +b^{T} \x +1 =\det\left(\bmatrix{I & C \x \\ \x^{T} C^{T} & 1+b^{T}\x}\right).
\end{equation}
Here the identity matrix $I$ is $r \times r$ matrix, where $r$ is the rank of $A$. Using Schur complement determinant formula for a partitioned matrix $L(\x):= \bmatrix{\mathbf{0} & C \x \\ (C \x)^{T} & b^{T} \x}$ the characteristic polynomial of $L(\x)$ is given by
\begin{eqnarray}
\det\left( \left[ \begin{array}{c:c}
    \lambda I_{r \times r} & C\x_{r \times 1} \\ 
    \hdashline
    (C\x)^{T}_{1 \times r} &  (\lambda - b^{T} \x)_{1 \times 1}
  \end{array} \right]\right) &=& \det(\lambda I) \det(\lambda-b^{T}\x-(C \x)^{T} (\lambda I)^{-1} (C \x)) \nonumber\\ &=& \lambda ^{r}(\lambda -b^{T}\x -\frac{1}{\lambda}\|C\x\|^{2}). \nonumber
\end{eqnarray}
Observe that at most two of $r+1$ eigenvalues are non-zero irrespective of the choice of $\x \in \R^{n}$. Further if there are exactly two non-zero eigenvalues (which implies rank($L(\x))=2$), these two non-zero eigenvalues are of opposite signs. This implies that there does not exist any $\x \in \R^{n}$ such that $x_{1}A_{1}+\dots+x_{n}A_{n} \succeq 0$ i.e., $L(\x) \nsucceq 0$ for any $\x \in \R^{n}$. Therefore, spectrahedron $S$ does not contain a full dimensional cone.

Observe that the above construction gives a symmetric MDR for a quadratic polynomial
$f(\x) = \x^{T}A\x + b^{T} \x + 1$, where the matrix $A$ is negative semidefinite. We can also get a
whole set of Hermitian MDRs from the above construction, by combining pairs of rows of the matrix $C$ where $A = -C^{T}C$. We illustrate this with an example.
\begin{example}
\rm{
Consider the polynomial $1-x_{1}^{2}-x_{2}^{2}-x_{3}^{2}-x_{4}^{2}-x_{5}^{2}$. Clearly, this polynomial
does not have an MDR of size $2$, since the rank of $Q/(1, 1)$ is $5$. On the other hand, using the
construction given above one obtains a size $6$ symmetric MDR given by the linear matrix
polynomial
\begin{equation*}
\bmatrix{1 & 0 & 0 & 0 & 0 & x_{1}\\0 & 1 & 0 & 0 & 0 & x_{2}\\0 & 0 & 1 & 0 & 0 & x_{3}\\0 & 0 & 0 & 1 & 0 & x_{4}\\0 & 0 & 0 & 0 & 1 & x_{5}\\x_{1}& x_{2}& x_{3}& x_{4}& x_{5} & 1}
\end{equation*}
One can now combine the first two rows of the $C$ matrix in this case and obtain a size $5$
MDR given by the linear matrix polynomial
\begin{equation*}
\bmatrix{1 & 0 & 0 & 0 & x_{1}+i x_{2}\\0 & 1 & 0 & 0 & x_{3}\\0 & 0 & 1 & 0 & x_{4}\\0 & 0 & 0 & 1 & x_{5}\\x_{1}-i x_{2}& x_{3}& x_{4}& x_{5} & 1}
\end{equation*}
Combining two other rows of the C matrix, one can go down to a size 4 Hermitian MDR.
For example,
\begin{equation*}
\bmatrix{1 & 0 & 0 & x_{1}+i x_{2}\\0 & 1 & 0 & x_{3}+ix_{4}\\0 & 0 & 1 & x_{5}\\x_{1}-i x_{2}& x_{3}-ix_{4}& x_{5} & 1}
\end{equation*}
Thus essentially using the same construction, one can build both symmetric and Hermitian
MDRs. Further observe that this process gives a whole range of sizes for the MDRs.
}
\end{example}
We now characterize all quadratic polynomials that exhibit an MDR.
\begin{theorem}\label{themcomcharquad}
A quadratic polynomial $f(\x)=\x^{T}A\x+b^{T}\x+1=Z^{T}[\x]QZ[\x] \in \R[\x]$ admits an MDR if and only if either one (possibly both) of the following two conditions is true.
\begin{enumerate}
\item $A$ is negative semidefinite
\item $Q/(1,1)$ is negative semidefinite and $\rank(Q/(1,1)) \leq 3$
\end{enumerate}
where $Q/(1,1)=A-\frac{1}{4}bb^{T}$. 
\end{theorem}
\pf If $A$ is negative semidefinite, then one can construct a symmetric MDR as demonstrated above to obtain  $\bmatrix{I & C\x \\ (C \x)^{T} & 1+b^{T}\x}$ where $A=-C^{T}C$. On the other hand if $Q/(1,1)$ is negative semidefinite and $\rank(Q/(1,1)) \leq 3$, then by the Theorem \ref{sufficientqmsdr} $f(\x)$ has an MDR of size $2$.

Conversely, if $f(\x)$ has an MDR, then $f(\x)$ is a RZ polynomial. By Proposition ~\ref{rznasc}, $Q/(1,1)=A-\frac{1}{4}bb^{T}$ is certainly negative semidefinite. So either $A$ is negative semidefinite or $A$ is not negative semidefinite.
If $A$ is not negative semidefinite, then by Lemma~\ref{indefinite} we know that the spectrahedron $S$ contains a full dimensional cone. So, in this case if the quadratic polynomial $f(\x)$ has an MDR of some size, then $f(\x)$ has an MDR of size $2$ by the Theorem \ref{fulldimensional}. On the other hand, by the Theorem \ref{sufficientqmsdr} if quadratic polynomial $f(\x)$ has an MDR of size $2$ then $Q/(1,1)$ is negative semidefinite and $\rank(Q/(1,1)) \leq 3$. Therefore, in the case of $A$ being not negative semidefinite, if $f(\x)$ has an MDR, then $Q/(1,1)$ is negative semidefinite and $\rank(Q/(1,1)) \leq 3$.  \hfill{$\square$}
\begin{remark}{\rm
The above theorems characterizes all quadratic polynomials that have an
MDR. The size of MDR for a quadratic polynomial $f(\x) = \x^{T}A\x + b^{T} \x + 1$ with $A$ being
negative semidefinite can range from $\lceil \frac{r}{2}+1 \rceil$ to $r + 1$, where rank $(A) = r$. Of course, even larger sizes MDRs
are possible, but these MDRs are such that the intersection of kernels of all the matrices Aj
would be non-trivial. Factoring out this common kernel, would reduce the situation to one
of the sizes outlined above. On the other hand, if the matrix $A$ is not negative semidefinite,
then $f(\x)$ has an MDR guarantees that $f(x)$ has an MDR of size $2$. This restricts the rank
of $A$ which is not negative semidefinite to a maximum of $4$, for an MDR to exist. In other
words, a polynomial $f(\x) = \x^{T}A\x + b^{T} \x + 1$ with rank$(A) > 4$, has an MDR if and only if
$A$ is negative semidefinite.}
\end{remark}
\begin{example}\rm{
We once again invoke Example \ref{exquniequiv}. Recall that condition $2$ of Theorem \ref{sufficientqmsdr}
is satisfied by all three examples, whereas condition $1$ of Theorem \ref{sufficientqmsdr} is only satisfied by the
first two cases. Thus for case $1$, the polynomial $1-8x_{1}x_{2}-4x_{1}x_{3}-100x_{2}^{2}-12x_{2}x_{3}-x_{3}^{2}-5x_{1}^{2}$ has a size $3$ symmetric MDR given by the linear matrix polynomial
\begin{equation*}
\bmatrix{1 & 0 & 11x_{1}/5+4x_{3}/5\\0 & 1 & 2x_{1}/5+10x_{2}+3x_{3}/5\\11x_{1}/5+4x_{3}/5 &2x_{1}/5+10x_{2}+3x_{3}/5 & 1}
\end{equation*}
Combining the two rows of the C matrix in this case, one can also obtain a size $2$ Hermitian
MDR given by
\begin{equation*}
\bmatrix{1 & \frac{11+2i}{5}x_{1}+10ix_{2}+\frac{4+3i}{5}x_{3}\\\frac{11-2i}{5}x_{1}-10ix_{2}+\frac{4-3i}{5}x_{3} & 1}
\end{equation*}
It is instructive to note that this is precisely one of the size $2$ linear matrix polynomials
obtained for case $1$, in the follow-up Example \ref{exquniequiv2}.
Similarly, for case $2$, the polynomial $1 + 4x_{1} + 10x_{2}-x_{1}^{2}- 2x_{1}x_{2}- x_{2}^{2}$ has another size $2$
MDR given by the linear matrix polynomial $\bmatrix{1 & x_{1}+x_{2}\\x_{1}+x_{2} &1+4x_{1}+10x_{2}}$ which is orthogonally
equivalent to MDR obtained in Example \ref{exquniequiv}.
}
\end{example}
As a result of the Theorem \ref{themcomcharquad} we have the following corollary.
\begin{corollary} \label{diagmsdrcor2}
A quadratic multivariate polynomial $f(\x)=\x^{T}A\x+b^{T}\x+1 \in \R[\x]$ has a diagonal MSDR of any size if and only if
\begin{enumerate}
\item Rank of $Q/(1,1) \leq 1$ and
\item $Q/(1,1)$ is a negative semidefinite matrix.
\end{enumerate}
\end{corollary}
\pf It follows from the equation (\ref{eqnsd}) that when $A$ is negative semidefinite, there can not exist any diagonal MSDR, otherwise $f(\x)$ can not be a quadratic polynomial. When $A$ is not negative semidefinite, but if quadratic polynomial has an MSDR of some size $>2$, it is proved in the Theorem \ref{themcomcharquad} that it has an MSDR of size $2$ too. It is shown in Corollary \ref{diagmsdrcor} that a quadratic polynomial has a diagonal MSDR of size $2$ if and only if $Q/(1,1)$ is a negative semidefinite matrix and of rank $\leq 1$. \hfill{$\square$}
\section{Conclusion}
It is well known that a quadratic polynomial $f(\x)=\x^{T}A\x+b^{T}\x+1$ has an MDR with $(n+1) \times (n+1)$ LMP, if the matrix $A$ is negative semidefinite. On the other hand, not much seems to be available in literature about quadratic polynomials having MDR of size $2$. In this chapter, we provide necessary and sufficient conditions for a quadratic multivariate polynomial to have an MDR of size $2$. We have provided an algorithm which can be used to construct MDRs of size $2$, when they exist. It has been shown that quadratic polynomials
having MDRs of size $2$ are of two kinds -- those that have exactly two unitarily non-equivalent Hermitian MDRs of size $2$ with none of the MDRs being a symmetric one and those that have exactly one unitarily equivalent Hermitian MDR which includes symmetric MDRs. It is further shown that all possible symmetric MDRs are orthogonally equivalent. Further,
we have shown that if a quadratic polynomial having $A$ which is not negative semidefinite has no MSDR of size $2$, it cannot have MSDR of any size greater than $2$. Consequently, we have completely characterized quadratic polynomials which have an MHDR (MSDR) of any size. The class of such quadratic polynomials belongs to any one of the following two categories: quadratic polynomials having $A$ which is negative semidefinite or quadratic polynomials which have an MDR of size $2$. As a consequence of this result, we have effectively characterized spectrahedra with non-empty interior defined by $2 \times 2$ linear matrix inequalities. 
Furthermore, if the mentioned conditions in Theorem \ref{themcomcharquad} are true, quadratic optimization problems (including (QCQP) and trust-region subproblems can be converted into an SDP relaxation problem irrespective of the fact that quadratic functions (objective or constraints) are convex or not.
\bibliographystyle{alpha}
\bibliography{refqplmi}
\end{document}